\documentclass[a4paper, 12pt]{amsart}
\usepackage{mathptmx, amssymb,amscd,latexsym, eulervm,calligra,mathrsfs,amsmath,amsthm,amsfonts,amsrefs,tikz,comment}
\usepackage[colorlinks=true,linkcolor=blue,urlcolor=blue,citecolor=blue]{hyperref}
\usepackage{xurl}
\hypersetup{breaklinks=true}
\usepackage[inner=2.8cm,outer=2.8cm, bottom=3.5cm]{geometry}
\allowdisplaybreaks

\newtheorem{thm}{Theorem}

\newtheorem{cor}[thm]{Corollary}

\theoremstyle{definition}

\newtheorem{problem}[thm]{Problem}
\newtheorem{question}[thm]{Question}

\newtheorem{remark}[thm]{Remark}

\numberwithin{equation}{section}

\newcommand{\Hilb}{\mathrm{Hilb}}
\newcommand{\HP}{\mathrm{HP}}
\newcommand{\HF}{\mathrm{HF}}
\newcommand{\Lex}{\mathrm{Lex}}

\renewcommand{\AA}{\mathbb{A}}
\renewcommand{\P}{\mathbb{P}}
\newcommand{\N}{\mathbb{N}}
\newcommand{\CC}{\mathbb{C}}
\newcommand{\kk}{\Bbbk}

\newcommand{\mfp}{{\mathfrak{p}}}

\newcommand{\mfh}{{\mathfrak{h}}}

\newcommand{\mcH}{{\mathcal{H}}}

\newcommand{\msfH}{{\mathsf{H}}}

\makeatletter
\@namedef{subjclassname@2020}{
  \textup{2020} Mathematics Subject Classification}
\makeatother

\begin{document}

\author[G. Farkas]{Gavril Farkas}
\address{Gavril Farkas: Institut f\"ur Mathematik, Humboldt-Universit\"at zu Berlin \hfill \newline
\indent Unter den Linden 6, 10099 Berlin, Germany}
\email{{\tt farkas@math.hu-berlin.de}}

\author[R. Pandharipande]{Rahul Pandharipande}
\address{Rahul Pandharipande: Department of Mathematics,
ETH Z\"urich \hfill \newline
\indent
R\"amistrasse 101, 8092 Z\"urich, Switzerland}
\email{{\tt rahul@math.ethz.ch}}

\author[A. Sammartano]{Alessio Sammartano}
\address{Alessio Sammartano: Dipartimento di Matematica, Politecnico di Milano \hfill\newline
\indent Via Bonardi 9, 20133 Milan, Italy}
\email{{\tt alessio.sammartano@polimi.it}}

\title{Irrational components of the Hilbert scheme of points}


\begin{abstract}
  We construct irrational irreducible components of the Hilbert scheme of points
  of affine space $\mathbb{A}^n$ for $n \geq 12$.
  We start with irrational components of the Hilbert scheme of curves in $\mathbb{CP}^3$ and use methods developed by Jelisiejew to relate these
  to irreducible components of the Hilbert schemes of points
of $\AA^n$. The result solves
  Problem XX of the Hilbert scheme of points problem list \cite{JJOpenProblems}.
\end{abstract}

\maketitle

\section{Introduction}

The Hilbert scheme $\Hilb (\AA^n)$ of points,
parametrizing finite subschemes of $\AA^n$,
plays an important role in current research in algebraic
geometry with connections to many fields including
commutative algebra, combinatorics, enumerative geometry, representation theory,
singularity theory, topology, and complexity theory.
The Hilbert scheme decomposes into connected components, $$\Hilb (\AA^n) = \bigsqcup_{d=0}^\infty \Hilb_d (\AA^n)\, ,$$
where each $ \Hilb_d (\AA^n)$ parametrizes finite schemes  of length $d$.
While $\Hilb_d(\AA^n)$ is smooth and irreducible for $n \leq 2$, the Hilbert
scheme of $d$ points is
 singular and reducible if $n \geq 3$ and $d$ is large.
For each $d$, there is a distinguished irreducible component of $ \Hilb_d (\AA^n)$, called the {\em smoothable} component,
whose general points parametrize reduced subschemes.
The smoothable component is generically
smooth of dimension $dn$ and is rational,
see \cite[Theorem C]{LellaRoggero}.
Very little is known about the geometry of the other components.

We investigate the classical problem of rationality.
Since components of $\Hilb (\AA^n)$ may  be non-reduced in their natural scheme structure,
any birational attribute of a component of $\Hilb (\AA^n)$ refers to the underlying reduced algebraic variety.
There are many results  known about the birational geometry of moduli spaces of  curves
 \cite{EisenbudHarris,HarrisMumford, FJP20}, $K3$ surfaces \cite{GHS07},  abelian varieties \cite{Mumford83}, and
other higher dimensional varieties.
The situation is markedly different for 0-dimensional objects,
as  their lack of intrinsic geometric structure  makes them harder to describe.
The birational type of non-smoothable components is understood in only a handful of specific examples.
In fact, all the components that have been described explicitly 
are obtained by constructing loci  $\mathcal{L}\subseteq \Hilb (\AA^n)$
which are isomorphic to open subsets of  products of  Grassmannians \cite[Remark 6.10]{JJElementary}.
In particular, all such constructions yield  rational components.
Moreover,
methods for certifying rationality based on Gr\"obner strata
\cite{LellaRoggero}
are successful in the case of Hilbert schemes of varieties of positive dimension,
but in the case of points they do not provide
information about components other than the smoothable one.

The following problem has been open for a long time. See, for example,
\cite[Problem XX]{JJOpenProblems} or
\cite[Problem 1.4]{AIM}.

\begin{problem}\label{ProblemRational}
  {\em Does there exist an irrational irreducible component  of the Hilbert scheme
    of points $\Hilb (\AA^n)$?}
\end{problem}

Our main result is a positive answer to
Problem \ref{ProblemRational}.
Specifically, over the field $\mathbb{C}$,
for each $n \geq 12$, we show the existence of
components of $\Hilb (\AA^{n})$
that are not rational (and not even rationally connected).
In fact,
we prove a stronger result: for each $n \geq 12$, the Hilbert scheme $\Hilb (\AA^{n})$
contains irreducible components with MRC-fibrations
of arbitrarily large dimension. The dimension of the MRC-fibration 
 can be understood as a measure of
irrationality, see Section \ref{prelims} for the definitions.
Our construction is inspired by  ideas related to the study of
Murphy's Law for singularities \cite{Vakil,Erman,Jelisiejew},
where pathologies are transferred across moduli spaces of different types.

In Section \ref{SectionReductionArtinian},
we perform a reduction from arbitrary graded Hilbert schemes{\footnote{By {\em graded Hilbert
      scheme}, we mean the Hilbert scheme of homogeneous ideals of a polynomial
    ring with respect to the standard degree \cite{HaimanSturmfels}. The definitions can be found in Section \ref{prelims}. }}
of a polynomial ring to
graded Hilbert schemes parametrizing finite graded algebras.
In Section \ref{SectionDominantMorphisms},
we construct dominant morphisms
from components of  $\Hilb (\AA^{n\geq 12})$
to the moduli space of curves in several steps starting with 
Hilbert schemes of curves in $\mathbb{CP}^3$ and
making essential use of Jelisiejew's TNT frames \cite{Jelisiejew}.

The bound $n\geq 12$ arises for us as follows.
We start with an irrational component of the Hilbert scheme of
$\mathbb{CP}^m$ for some $m$. All such Hilbert schemes are rational
for $m\leq 2$, so  $\mathbb{CP}^3$   is the smallest dimensional
projective space that we can use
(with irrational components
obtained from Hilbert schemes of curves). A non-empty open set of
a component of the Hilbert scheme
of curves in $\mathbb{CP}^3$ is an open set of a graded Hilbert scheme in 4 homogeneous
variables $x_0,x_1,x_2, x_3$. In order to use Jelisiejew's TNT frames
to show that the latter is an open set of a component of $\Hilb (\AA^{n})$,
we must first add two more variables $x_4$ and $x_5$ and then
double the total number of variables. So, $12$ arises as $2\cdot (4+2)$. 
Whether irrational components occur
in $\Hilb (\AA^{n<12})$ is an interesting question.

Another question concerns the number of points $d$ required for 
$\Hilb_d (\AA^{12})$ to have an irrational component.
The genus and the degree of the curves parametrized by the irrational
component of the Hilbert scheme of $\mathbb{CP}^3$ determine
a corresponding point number $d$ for $\AA^{12}$. 
We have made no
effort to optimize (or even to bound) the very large $d$
required for our construction.

\section{Preliminaries}
\label{prelims}

Let $\kk$ denote an arbitrary  field. 
Let $S= \kk[x_0, \ldots, x_m]$ be the
polynomial ring in $m+1$ variables.
If $M$ is a finitely generated graded $S$-module, the Hilbert function
$$\HF(M)\colon \N\to \N$$
is defined by $\HF(M;d) = \dim_\kk [M]_d$,
where $[M]_d$ denotes the graded component of degree $d$.
The Hilbert function agrees with the
Hilbert polynomial $\HP(M;z) \in \mathbb{Q}[z]$ for sufficiently large $d\in\N$.

Let $m \in \N$ and $\mfp(z) \in \mathbb{Q}[z]$.
We denote by $\Hilb_{\mfp(z)}(\P^m)$
the Hilbert scheme parametrizing closed subschemes of $\P^m$ with Hilbert polynomial $\mfp(z)$.
Equivalently, $\Hilb_{\mfp(z)}(\P^m)$
parametrizes saturated homogeneous ideals $$I\subseteq S= \kk[x_0, \ldots, x_m]$$
such that the algebra $S/I$ has Hilbert polynomial $\mfp(z)$.
For simplicity,
we use the same symbol $I$ to denote both the ideal and the associated
$\kk$-point of the Hilbert scheme,
$$I \in \Hilb_{\mfp(z)}(\P^m)\, .$$

Let $S$ be a polynomial ring over $\kk$, and let $\mfh \colon \N \to \N$  be a  function.
We denote by $\mcH^\mfh(S)$ the graded Hilbert scheme parametrizing
homogeneous ideals $I\subseteq S$ such that the algebra $S/I$ has
Hilbert function equal to $\mfh$ \cite{HaimanSturmfels}.

A projective variety  $X$ over $\CC$ is {\em rationally connected}
if, through a pair $x, y$ of general points of $X$, there exists a rational curve $f\colon \mathbb P^1\rightarrow X$ passing through $x$ and $y$. Unirational varieties are rationally connected and rational connectedness is clearly a birational invariant.
We refer to \cite[Section IV.3]{Kollar} for equivalent definitions  and related notions.

Let $X$ be a smooth proper variety over $\CC$.
The  maximal rationally connected fibration, or \emph{MRC-fibration} \cite{campana}, \cite[Definition IV.5.3]{Kollar}, of $X$
is a dominant rational map $$\varphi\colon  X \dashrightarrow MRC(X)$$
with  rationally connected fibers
which is maximal
with respect to the property of having rationally connected fibers. Equivalently, for a very general point $x\in MRC(X)$ any rational curve in $X$ passing through a general point of $\varphi^{-1}(x)$ is contained in the fibre $\varphi^{-1}(x)$.
The variety $MRC(X)$ is called the \emph{MRC-quotient} of $X$.
The MRC-fibration  is uniquely determined up to birational equivalence
and is functorial for dominant maps.
The variety $MRC(X)$ is a point if and only if $X$ is rationally connected, and
$MRC(X)=X$ if and only if $X$ is not uniruled.
The MRC dimension of $X$ is defined as  $\dim (MRC(X))$ and is
clearly a birational invariant of $X$.
If $X$ is not smooth or proper,
we  extend the notion above by taking a desingularization of a compactification of $X$.
See \cite[Section IV.5]{Kollar} for further properties of MRC-fibrations.

\section{Reduction to finite graded algebras}\label{SectionReductionArtinian}

Let $\kk$ be an algebraically closed field.
We prove here a general result of independent interest:
every  graded Hilbert scheme  of a  polynomial ring over $\kk$  is isomorphic to a
graded Hilbert scheme that parametrizes finite graded algebras over $\kk$.

We begin by recalling some definitions and facts about Hilbert functions and  syzygies.
Let $S = \kk[x_0,\ldots, x_m]$ be a polynomial ring.
A monomial ideal $L\subseteq S$ is a {\em lexsegment ideal} if every graded component of $L$
is spanned by an initial segment of monomials of $S$, where monomials are ordered lexicographically.
By Macaulay's Theorem \cite[Theorem 4.2.10]{BrunsHerzog}, there is a bijection  between Hilbert functions of quotient algebras of $S$ and lexsegment ideals of $S$.
If $\mfh = \HF(S/I)$ for some homogeneous ideal $I\subseteq S$, we denote by $\Lex(\mfh)$ the unique lexsegment ideal of $S$ with $\mfh = \HF(S/\Lex(\mfh))$.

Let $I\subseteq S$ be a homogeneous ideal.
The integers  $\beta_{i,j}(I) := \dim_\kk [\mathrm{Tor}_i(I,\kk)]_j$
are the graded Betti numbers of $I$.
The Castelnuovo-Mumford regularity is defined as the quantity
$$\mathrm{reg}(I) = \max \{\, j-i \,\mid\, \beta_{i,j}(I) \ne 0\, \}\, .$$
The degrees of a minimal set of generators of $I$ are bounded by $\mathrm{reg}(I)$.
We have $$\HF(S/I;d) = \HP(S/I;d)$$ for all $d \geq \mathrm{reg}(I)$.
By the Bigatti-Hulett-Pardue Theorem \cite[Theorem 31]{Pardue},
the lexsegment ideal $\Lex(\mfh)$ attains the largest Betti numbers $\beta_{i,j}$ among all homogeneous ideals with Hilbert function $\mfh$, for every $i,j$.

For an integer $B \in \N$, we denote by $I_{\geq B}$, respectively, by  $I_{\leq B}$,
the ideal generated by all the homogeneous polynomials of $I$ of degree at least $B$,
respectively,
 at most $B$.

Our main result concerning the reduction to finite graded algebras is the following.

\begin{thm}\label{TheoremArtinian}
Let $S = \kk[x_0,\ldots, x_m]$ and let $\mfh\colon \N\to \N$ be a  function.
For $D \in \N$, let $\overline{\mfh}:\N\to \N$ be the function defined by
$$
\overline{\mfh}(d) = \begin{cases}
\mfh(d) & \text{ if } d < D, \\
0 & \text{ if } d \geq D.
\end{cases}
$$
Then, for all
 $D \gg 0$,  there is  an isomorphism between the graded Hilbert schemes
$\mcH^\mfh(S)$
and
$\mcH^{\overline{\mfh}}(S)$.
\end{thm}

\begin{proof}
We begin by constructing a morphism $\theta :
\mcH^\mfh(S) \to \mcH^{\overline{\mfh}}(S)$.

Let $K$ be an arbitrary field.
By the Bigatti-Hulett-Pardue Theorem \cite[Theorem 31]{Pardue},
there exists an integer $D_0 \in \N$, depending only on $m$ and $\mfh$,
such that $\mathrm{reg}(I) \leq D_0$ for all homogeneous ideals $I \subseteq K[x_0,\ldots, x_m]$ whose quotient algebra has Hilbert function $\mfh$.

By  \cite[Proposition 3.8]{BayerMumford} or
\cite[Corollary 2.7]{CavigliaSbarra},
there exists an integer $D_1\in \N$, depending only on $D_0$ and $m$,
such that we have $\mathrm{reg}(I) \leq D_1$ 
for all homogeneous ideals $I \subseteq K[x_0,\ldots, x_m]$
generated by polynomials of degree at most $D_0$.

Choose an integer $D \geq D_1+m+1$, and let $\overline{\mfh}$ be the truncation of $\mfh$ in degree $D$, defined as in the statement of the theorem.
The Hilbert scheme $\mcH^\mfh(S)$ 
represents the Hilbert functor $\underline{\mathrm{Hilb}}^\mfh_S$,
from the category of $\kk$-algebras to the category of sets, defined by 
$$
R \mapsto \Big\{\text{homogeneous } I \subseteq R\otimes_\kk S \, \big| \,  [(R\otimes_\kk S)/I]_d \text{ locally free over } R \text{ of rank } {\mfh}(d)\,\, \forall \,\, d \Big\},
$$
and likewise for  $\overline{\mfh}$,
see \cite{HaimanSturmfels} for details.
Let  $\mathbf{m} = (x_0, \ldots, x_m)$ denote the homogeneous maximal ideal of $R \otimes_\kk S$.
For a given $I \in \underline{\mathrm{Hilb}}^\mfh_S(R)$,  the formula 
\begin{equation}\label{EqThetaMap}
\theta(I) = I + \mathbf{m}^D 
\end{equation}
yields an ideal
 $\theta(I) \subseteq R\otimes_\kk S$ such that $ [(R\otimes_\kk S)/\theta(I)]_d= [(R\otimes_\kk S)/I]_d$ 
 if $d < D$ and
$ [(R\otimes_\kk S)/\theta(I)]_d=0$ otherwise.
In particular,  $ [(R\otimes_\kk S)/\theta(I)]_d$ is locally free of rank $\overline{\mfh}(d)$ for all $d$.
Thus, $\theta$ defines a natural transformation $ \underline{\mathrm{Hilb}}^\mfh_S(R) \to  \underline{\mathrm{Hilb}}^{\overline{\mfh}}_S(R)$, 
and a morphism  of schemes $\theta: \mcH^\mfh(S) \to \mcH^{\overline{\mfh}}(S)$.

We  show that $\theta$ is bijective on $K$-points, where $\kk \subseteq K$ is any field extension.
Observe that a $K$-point of $\mcH^\mfh(S)$ is a homogeneous ideal 
$$
I \subseteq K \otimes_\kk S = K[x_0,\ldots,x_m] =:T
$$ 
whose quotient has Hilbert function $\mfh$,
that is, an ideal $I\in \mcH^\mfh(T)$;
likewise for a $K$-point of $\mcH^{\overline{\mfh}}(S)$.
We use the symbol $\theta$ also to denote the induced map $\theta: \mcH^\mfh(T) \to \mcH^{\overline{\mfh}}(T)$,
which is defined by the same formula \eqref{EqThetaMap}.
Consider the map  $\rho \colon \mcH^{\overline{\mfh}}(T) \to  \mcH^\mfh(T)$ defined, 
at the level of $K$-points, by the formula 
$
\rho(J) = J_{\leq D_0}.
$
We claim that $\rho$ and $\theta$ are inverses as set-theoretic maps.

For every ideal $I \in \mcH^\mfh(T)$,
 since $I$ is generated in degrees at most equal to $D_0$ and  since $D_0<D$, we have
$$I = I_{\leq D_0} = (I+\mathbf{m}^{D})_{\leq D_0}\, ,$$ 
therefore, $\rho(\theta(I)) = I$.
Consider, in particular, the
lexsegment ideal with Hilbert function $\mfh$
$$L =\Lex(\mfh) \in \mcH^\mfh(T).$$ 
Then, 
$$\overline{L}= \theta(L) = L+\mathbf{m}^{D}\in \mcH^{\overline{\mfh}}(T)$$
is again a lexsegment ideal, with Hilbert function equal to  $\overline{\mfh}$.
Since $L$ is generated in degrees up to $D_0$,  $\overline{L}$ is generated in degrees up to $D_0$ and in degree $D$.
Using the Bigatti-Hulett-Pardue Theorem again, we conclude that 
every ideal $J \subseteq T$ with Hilbert function equal to ${\overline{\mfh}}$ is generated, at most, in degrees up to $D_0$ and in degree $D$.

Consider now any such ideal $J \subseteq T$ with Hilbert function equal to ${\overline{\mfh}}$.
By the conclusion of the previous paragraph, the graded components of the ideals $J$ and $\rho(J) =J_{\leq D_0}$ agree up to degree $D-1$.
In other words, $$\HF(T/\rho(J);d) =  \overline{\mfh} (d)= \mfh (d) $$
for all $d\leq D-1$
(in particular, for all $d \leq D_1+m$).
Since $\rho(J)$ is generated in degrees up to $D_0$, it follows,
by choice of $D_1$,
 that the regularity of $\rho(J)$ is at most equal to $D_1$.
Therefore,
$$\HP(T/\rho(J);d) = \HF(T/\rho(J);d)$$ for all $d \geq D_1$.
Let $\mfp(z)$ be the Hilbert polynomial of the Hilbert function $\mfh$.
By choice of $D_0$, we have $\mfp(d) = \mfh(d)$ for all $d \geq D_0$.
In conclusion, we have
$$
\HP(T/\rho(J);d) = \HF(T/\rho(J);d) =   \overline{\mfh} (d)= \mfh (d) = \mfp(d)
$$
for all $d = D_1, D_1+1, \ldots, D_1+m$.
Since $\HP(T/\rho(J);d) $ and $ \mfp(d)$ are univariate polynomials of degree at most $m$,
they must be the same polynomial.
We deduce that
$$\HF(T/\rho(J);d) =\mfp(d) = \mfh(d) $$ for all $d \geq D_1$, in other words,
 $T/\rho(J)$ has Hilbert function equal to $ \mfh $, and $\rho(J) \in \mcH^{{\mfh}}(T)$.
 
  Since $J$ is generated in degrees up to $D_0$ and in degree $D$ and $J_{\geq D} = \mathbf{m}^{D}$, we obtain
 $$J = J_{\leq D_0} + \mathbf{m}^{D}\, ,$$ 
 therefore, $\theta(\rho(J))=J$ as desired.
Therefore, we have proven that
$\theta: \mcH^\mfh(S) \to \mcH^{\overline{\mfh}}(S)$ is  bijective on $K$-points.

Next, we show that $\theta$ is an open immersion.
By \cite[Proposition 3.1]{Erman}, for each $I \in \mcH^\mfh(S)$,
the morphism $\theta$ induces an isomorphism between the functors of Artin rings 
which send a local Artinian $\kk$-algebra $R$ to the set of homogeneous deformations over $R$ of $S/I$ and of $S/\theta(I)$, respectively.
These two functors of Artin rings  are pro-represented by the completed local rings 
$\widehat{\mathscr{O}}_{\mcH^\mfh(S),I}$ and 
$\widehat{\mathscr{O}}_{\mcH^{\overline{\mfh}}(S),\theta(I)}$
of the two Hilbert schemes at the respective $\kk$-points (see for example \cite[Section 2.2]{Sernesi}).
Thus, $\theta$ induces isomorphisms of completed local rings at $\kk$-points
and therefore is an étale morphism (since $\kk$ is algebraically closed).

 Since $\theta$ is injective on $K$-points for any
field $\kk \subseteq K$,
it follows from \cite[Lemma 29.10.2]{stacks} that it is universally injective.
Since $\theta$ is universally injective and étale, 
it is an open immersion by  \cite[Lemma 35.25.2]{stacks}.

Finally, 
since
$\theta$ is an open immersion which is bijective on $\kk$-points,
and $\kk$ is algebraically closed, 
it follows that it must be an
isomorphism.
\end{proof}

\section{Dominant maps between components of Hilbert schemes}\label{SectionDominantMorphisms}

We prove here the main result of the paper
by constructing dominant rational maps between irreducible components  of various Hilbert schemes.
An important tool for us is   the theory of TNT frames of \cite{Jelisiejew},
which produces dominant rational maps in the form of local retractions.
A {\em retraction} is a morphism of schemes $\pi\colon X\to Y$ with a section,
$$\iota \colon Y\to X \ \ \ \text{such that} \ \ \
\pi \circ \iota = \mathrm{id}_Y\, .$$
A {\em local retraction} of pointed schemes $(X,x) \to  (Y,y)$ is a retraction $(U,x) \to (V,y)$ for open subsets
$x \in U \subseteq X$ and $y \in V \subseteq Y$.

\vskip 4pt

\begin{thm}\label{MainTheorem}
  There exist irreducible components of $\Hilb(\AA^{12})$ over $\CC$
  which are not rationally connected.
Moreover,  there exist irreducible components of $\Hilb(\AA^{12})$ of arbitrarily large MRC dimension.
\end{thm}

\begin{proof}
Let $\mathscr{M}_{g}$ be the moduli space of curves of genus $g$.
By \cite{EisenbudHarris,HarrisMumford, FJP20},
 $\mathscr{M}_{g}$ is of general type if $g \geq 22$.
We will prove the theorem by constructing a dominant rational  map
$$\mathcal{C}\dashrightarrow \mathscr{M}_{g}\, ,$$
where $\mathcal{C}\subseteq \Hilb(\AA^{12})$ is an irreducible component.
We will realize the  map
as a composition of a series of dominant  morphisms $(\msfH_{i+1},I_{i+1})
\to
(\msfH_{i},I_{i})$,
where each $\msfH_i$ is an irreducible open subscheme of a corresponding
Hilbert scheme and $I_i\in \msfH_i$ is an ideal parametrized by that Hilbert scheme.

Fix integers $g \geq 22$ and $d\geq \frac{3g+12}{4}$, 
and consider smooth curves in $\P^3$ of genus $g$ and degree $d$.
The inequality $d\geq \frac{3g+12}{4}$ amounts to the non-negativity of the Brill--Noether number $\rho(g, 3,d)=g-4(g-d+3)
$.
 It follows from the 
main results of Brill--Noether theory
that there exists a unique irreducible component 
$$Y \subseteq \Hilb_{dz+1-g}(\P^3)$$
which dominates the moduli space $\mathscr{M}_{g}$,  
see e.g. \cite{EisenbudHarris2} and references therein or
 \cite[Theorem 19.9]{EisenbudHarrisBook}.
For a general curve $C\subseteq \mathbb P^3$ parametrized by $Y$, the Hilbert function of the coordinate ring
is determined by the \emph{Maximal Rank Theorem} \cite{BE,Larson} (see
also \cite[Theorem 12.10]{EisenbudHarrisBook}), asserting that  multiplication maps
$$\mbox{Sym}^a H^0(C, \mathscr{O}_C(1))\rightarrow H^0(C, \mathscr{O}_C(a))$$ are  of maximal rank for all integers $a\geq 1$. 
In particular, there exists a non-empty open subset $\msfH_1 \subseteq Y$ where the Hilbert function of $S/I$ is constant for all $I \in \msfH_1$.
By comparing dimensions,
 curves $C$ parametrized by $\msfH_1$ {\em never} lie on a  quadric surface.
Let $\mfh$ denote this Hilbert function, and let $I_1 \in \msfH_1$ be an ideal that belongs to no other component of $\Hilb_{dz+1-g}(\P^3)$.

Let $S= \CC[x_0, x_1,x_2, x_3]$ and
consider the saturated ideal $I_2 = I_1$ as a point in the graded Hilbert scheme $\mcH^\mfh(S)$.
The natural map $\mcH^\mfh(S) \to \Hilb_{dz+1-g}(\P^3)$ restricts to an isomorphism $(\msfH_2,I_2) \cong (\msfH_1,I_1)$
for some $\msfH_2\subseteq \mcH^\mfh(S)$ open and irreducible.

Let $P = S[x_{4}, x_{5}] = \CC[x_0, \ldots, x_{5}]$,
let
 $I_3 = I_2 \cdot P$ be the extended ideal,
and let $\mfh'  =\HF(P/I_3)$.
We claim that 
there is a local retraction
\begin{equation} \label{retrct}
  (\mcH^{\mfh'}(P),I_3) \to  (\mcH^{\mfh}(S),I_2)\, .
  \end{equation}

The local retraction \eqref{retrct} can be constructed by means of a 
suitable Bia\l{}ynicki--Birula decomposition, 
 proceeding  very similarly  to \cite[Proof of Theorem 1.3]{Jelisiejew}.
Specifically, let 
$\mathbb{G}_m = \CC^*$ act on $P$ with weight 1 on the variables $x_4, x_5$ and weight 0 on the remaining variables.
Let $\delta_1$ denote the standard grading $\delta_1(x_i) = 1$ for all $i$,
and let $\delta_2$ denote the $\mathbb{Z}$-grading on $P$ induced by the $\mathbb{G}_m$-action,
$$\delta_2(x_0)= \delta_2(x_1)= \delta_2(x_2)=\delta_2(x_3)= 0\, , \ \ \ \ 
\delta_2(x_4)=\delta_2(x_5) =1\, .$$
By general theory \cite{Drinfeld,JJElementary,JS},
there exists a $\CC$-scheme $\mcH^{\mfh'}(P)^-$,
called the negative{\footnote{The negative decomposition is the Bia\l{}ynicki--Birula decomposition associated to
the compactification $\overline{\mathbb{G}}_m = \mathbb{G}_m\cup \{\infty\}$,
that is, associated to the limit $t \to \infty$. 
It is the same choice  used in \cite{JJElementary,Jelisiejew}, however, it is there denoted by 
a $+$ sign instead of a $-$ sign. }}
Bia\l{}ynicki--Birula decomposition for 
$\mcH^{\mfh'}(P)$.
The scheme $\mcH^{\mfh'}(P)^-$
represents the functor of $\mathbb{G}_m$-equivariant families that admit a limit at $\infty$:
the functor assigns to a $\CC$-algebra $R$ the  ideals of $R [t^{-1},x_0, \ldots, x_{5}]$ that are homogeneous
with respect to $\delta_1$, with locally free components of  ranks $\mfh'(d)$, 
and are also homogeneous  with respect to $\delta_2$, where
$$\delta_1(t^{-1})=0 \ \text{  and  }\ \delta_2(t^{-1})=-1\, .$$
There is a natural map $\mcH^{\mfh'}(P)^-\to \mcH^{\mfh'}(P)$ by evaluating at $t=1$.
By  taking the limit $t \to \infty$, there is a natural map to the $\mathbb{G}_m$-fixed locus,
\begin{equation} \mcH^{\mfh'}(P)^-\to \mcH^{\mfh'}(P)^{\mathbb{G}_m}\, .  \label{mtofix}\end{equation}
Since there is a section  $\mcH^{\mfh'}(P)^{\mathbb{G}_m} \to \mcH^{\mfh'}(P)^- $ defined by taking constant families,
the map \eqref{mtofix}
is a retraction.
See \cite[Diagram (1.1)]{Jelisiejew} and the related discussion for more details. 

Since $I_3$ is homogeneous with respect to both $\delta_1$ and $\delta_2$, we can consider $I_3$
as a $\CC$-point of all three schemes  $\mcH^{\mfh'}(P),$ $\mcH^{\mfh'}(P)^-$, and $\mcH^{\mfh'}(P)^{\mathbb{G}_m}$.
The tangent space to $ \mcH^{\mfh'}(P)$ at the point $I_3$ is the vector space $[\mathrm{Hom}_P(I_3,P/I_3)]_0$ of
 $P$-linear homomorphisms that are of $\delta_1$-degree 0.
 The tangent space to  $ \mcH^{\mfh'}(P)^-$ 
 at  $I_3$
 is the subspace consisting of homomorphisms that are also
$\delta_2$-graded, of non-negative degree.
Since the generators of  $I_3$ have $\delta_2$-degree 0, 
and $P/I_3$ is concentrated in non-negative $\delta_2$-degrees,  
the two tangent spaces coincide, 
that is, $[\mathrm{Hom}_P(I_3,P/I_3)]_0$
is already concentrated in non-negative $\delta_2$-degrees.
As a consequence, the map
\begin{equation} \label{opim}
  \mcH^{\mfh'}(P)^-\to \mcH^{\mfh'}(P)
  \end{equation}
is an open immersion around $I_3$,
see 
\cite[Example 1.3 and Proposition 1.6]{JS} or  the argument of \cite[Theorem 4.5]{JJElementary}.
After combining the open immersion \eqref{opim} with the retraction \eqref{mtofix},
we obtain a local retraction
$$(\mcH^{\mfh'}(P),I_3)\to (\mcH^{\mfh'}(P)^{\mathbb{G}_m},I_3)\,.$$

For the last step in the construction of \eqref{retrct}, we
observe that $(\mcH^{\mfh'}(P)^{\mathbb{G}_m},I_3)$ is isomorphic to 
$(\mcH^{\mfh}(S),I_2)$,
with isomorphisms given by extension and contraction of ideals along the ring map $S \to P$.
The ideals of $P$ that are homogeneous with respect to both $\delta_1$ and $\delta_2$ and have the same bigraded Hilbert function as $I_3 = I_2 \cdot P$ are precisely the extensions of homogeneous ideals of $S$ with
the same Hilbert function as $I_2$.

Using the local retraction \eqref{retrct}, we conclude that  there is a  dominant  morphism
$$(\msfH_3, I_3) \to (\msfH_2,I_2)\, , $$
for some open irreducible $\msfH_3 \subseteq \mcH^{\mfh'}(P)$.
Since $I_2$ is saturated, we have $\mathrm{depth}(S/I_2) \geq 1$ and, therefore, $\mathrm{depth}(P/I_3) \geq 3 $.
Moreover, $I_3$ contains no quadrics by construction.

Let $D\in \N$ be an integer obtained by applying Theorem \ref{TheoremArtinian} to the  graded Hilbert scheme $\mcH^{\mfh'}(P)$,
and let
$$I_4 = I_3 + (x_0, \ldots, x_{5})^{D} \subseteq P\, .$$
By Theorem \ref{TheoremArtinian},
there exists an isomorphism $(\msfH_4,I_4) \cong (\msfH_3,I_3)$
for some open irreducible $\msfH_4 \subseteq \mcH^{\overline{\mfh}'}(P)$,
where $\overline{\mfh}' = \HF(P/I_4)$.

Let $T = P[y_0, \ldots, y_{5}] = \CC[x_0, \ldots, x_{5 }, y_0,\ldots, y_{5}]$,
and consider the ideal 
$$
I_5 = I_3 T + (x_0, \ldots, x_{5})^D+ (y_0, \ldots, y_{5})^2 + \left( \sum_{i=0}^5 x_i y_i\right) \subseteq T.
$$
Then, $T/I_5$ is a finite  $\CC$-algebra, so
$I_5 \in \Hilb(\AA^{12}).$
In \cite[Equation (1.3)]{Jelisiejew},
the ideal $I_5$ is called \emph{a TNT frame of size $a=D-1$ for $I_3$}.
Since  $\mathrm{depth}(P/I_3) \geq 3 $ and $I_3$ contains no quadrics,
we can apply \cite[Proposition 1.4]{Jelisiejew} and obtain
 a local retraction
 $$(\Hilb (\AA^{12}), I_5) \to  (\mcH^{\overline{\mfh}'}(P), I_4)\, .$$
Therefore,
there exists a dominant morphism $(\msfH_5,I_5) \to (\msfH_4,I_4)$ for some open irreducible
$ \msfH_5 \subseteq \Hilb (\AA^{12})$.

Finally, consider the irreducible component
$\mathcal{C}= \overline{\msfH}_5 \subseteq \Hilb (\AA^{12})$.
The composition of  the   maps
$$
(\msfH_5, I_5) \to (\msfH_4,I_4)\cong (\msfH_3,I_3) \to (\msfH_2,I_2) \cong (\msfH_1,I_1)
$$
together with the dominant rational map
$Y \dashrightarrow \mathscr{M}_{g}$ yields a dominant rational map
$$\mathcal{C}\dashrightarrow \mathscr{M}_{g}\, .$$
We obtain a dominant rational map  $\mathcal{C}_{red}\dashrightarrow \mathscr{M}_{g}$ on the reduced structure.
Since the moduli space $\mathscr{M}_{g}$ is of general type for $g \geq 22$, we conclude that $\mathcal{C}_{red}$ is not rationally connected.
In fact,
by  functoriality of MRC-fibrations \cite[Theorem IV.5.5]{Kollar},
we obtain a dominant rational map
$$MRC(\mathcal{C}_{red})\dashrightarrow MRC(\mathscr{M}_{g})\, .$$
Since $\mathscr{M}_{g}$ is of general type, we have
$MRC(\mathscr{M}_{g})= \mathscr{M}_{g}$.
Thus, the MRC-dimension of $\mathcal{C}_{red}$ is at least $\dim \mathscr{M}_{g} = 3g-3$.
\end{proof}

\begin{remark}
  While much of the argument of the proof of Theorem \ref{MainTheorem} is valid in arbitrary characteristic,
  the conclusion follows only in characteristic 0, 
since the geometry of moduli spaces is less understood in positive characteristic,
and, more generally, birational geometry becomes considerably more pathological.
\end{remark}

\begin{cor}\label{CorollaryHigher}
  The conclusion of Theorem
\ref{MainTheorem}
holds for $\Hilb(\AA^{n})$ for all $n \geq 12$.
\end{cor}

\begin{proof}
We continue to use the notation of the proof of Theorem \ref{MainTheorem}.
The argument is, again, based on a suitable 
Bia\l{}ynicki--Birula decomposition, and is very similar to the one used in the proof of Theorem \ref{MainTheorem}.
Specifically,
let $Q = T[z_{13}, \ldots, z_n]$ and consider the ideal
$$I_6 = I_5 \cdot Q+ (z_{13}, \ldots, z_n) \subseteq Q\, .$$
Let $\mathbb{G}_m = \CC^*$ act on $Q$ with weight 1 on the variables $z_i$ and weight 0 on the remaining variables.
The action induces a $\mathbb{Z}$-grading $\delta_3$ on $Q$ with 
$$\delta_3(x_i) = \delta_3(y_i) =0\, ,\ \  \delta_3(z_j) =1\, .$$
We consider now the positive 
Bia\l{}ynicki--Birula decomposition 
$\Hilb(\AA^n)^+$ for 
$\Hilb(\AA^n)$,
that is, 
the decomposition associated with
the compactification $\overline{\mathbb{G}}_m = \mathbb{G}_m\cup \{0\}$ 
associated to the limit $t \to 0$.

Since $I_6$ is generated in $\delta_3$-degrees 0 and 1 and $Q/I_6$ is concentrated in $\delta_3$-degree 0,
the tangent space $\mathrm{Hom}(I_6,Q/I_6)$ to $\mathrm{Hilb}(\AA^n)$ at the point $I_6$  
is concentrated in non-positive $\delta_3$-degrees.
Therefore, the tangent spaces to $\mathrm{Hilb}(\AA^n)^+$ and to $\mathrm{Hilb}(\AA^n)$ at the point $I_6$ coincide.
As before, it follows by \cite[Proposition 1.6]{JS} that the natural map
$$\mathrm{Hilb}(\AA^n)^+ \to \mathrm{Hilb}(\AA^n)$$
is an open immersion near $I_6$,
and that we have a local retraction to the fixed locus
$$(\mathrm{Hilb}(\AA^n),I_6) \to (\mathrm{Hilb}(\AA^n)^{\mathbb{G}_m},I_6)\, .$$
The connected component of $\mathrm{Hilb}(\AA^n)^{\mathbb{G}_m}$ containing $I_6$ parametrizes $\delta_3$-homogeneous ideals of $Q$ with the same Hilbert function as $I_6$.
Since they define quotients of $Q$ that are concentrated in $\delta_3$-degree 0,  they are exactly the ideals of the form $$I \cdot Q + (z_{13}, \ldots, z_n) \subseteq Q$$
where $I \in \Hilb(\AA^{12})$ has the same colength as $I_5$.
The connected component of $\mathrm{Hilb}(\AA^n)^{\mathbb{G}_m}$ containing $I_6$
is therefore isomorphic to the connected component of $\mathrm{Hilb}(\AA^{12})$ containing $I_5$.
We thus obtain a local retraction
$$(\mathrm{Hilb}(\AA^n),I_6) \to (\mathrm{Hilb}(\AA^{12}),I_5)\, ,$$
and hence
a dominant morphism $(\msfH_6,I_6) \to (\msfH_5,I_5)$ for some open irreducible
$ \msfH_6 \subseteq \Hilb (\AA^{n})$.
The conclusion follows as in the last paragraph of the proof of Theorem \ref{MainTheorem}
 by taking the composition with the dominant rational map  $\mathcal{C}\dashrightarrow \mathscr{M}_{g}$.
\end{proof}

\section{Questions}

It is unclear to us
whether the components $\mathcal{C}$ produced in Theorem \ref{MainTheorem} are
reduced.
A natural approach, by \cite[Theorem 4.6]{JJElementary},
 would involve showing the vanishing of the obstruction space $T^2(T/I_5)_{\geq 0}$ to prove the smoothness of the point $I_5 \in \mathcal{C}$.
However,
controlling the obstruction space is in general difficult.
See, for example, the delicate arguments of \cite{StaalSatriano}
for  much simpler ideals.

\begin{question}
  Are any of the irreducible components of $\mathrm{Hilb}(\AA^{12})$
  constructed in
  the proof of Theorem \ref{MainTheorem} reduced?
\end{question}

Though dimension 12 in Theorem \ref{MainTheorem}
is required for our method of proof, we would expect that 
$\mathrm{Hilb}(\AA^{n})$ contains irrational components also for  smaller $n$.
It would be especially interesting to see if irrational components exist already in $\mathrm{Hilb}(\AA^{3})$: dimension 3 is not only
the smallest dimension where the question is open, it is an
exceptional boundary case for Hilbert schemes in many respects
\cite{BBS,BF,GGGL,JKS,JRS,MNOP,RS,Ricolfi}.
Thus, we ask the following refinement of Problem \ref{ProblemRational}.

\begin{question}
Does $\mathrm{Hilb}(\AA^{3})$ contain irrational components?
\end{question}

Following \cite{BDELU17}, the degree of irrationality $\mbox{irr}(X)$ of an $n$-dimensional projective variety $X$ is the minimal degree of a dominant map $X\dashrightarrow \mathbb P^n$. We can ask the following question:

\begin{question}
Does $\mathrm{Hilb}(\AA^n)$ contain components of unbounded degree of irrationality?
\end{question}

The answer to the parallel question for the moduli space $\mathscr{M}_{g}$ is not known (and
expected to be out of reach at present), see also \cite[Problem 4.4]{BDELU17}. It is conceivable
that the case of $\mathrm{Hilb}(\AA^n)$ is more approachable.

In a more speculative direction,
we can ask  about the array of different birational types found in Hilbert schemes of points.
It is clear that not all birational types can occur:
for example,  components contain many rational curves and are never of general type.
A parallel situation holds for Murphy's law for singularities \cite{Vakil}:
we cannot expect to find all singularities on a moduli space, but rather, all
singularity types (equivalence classes induced by smooth morphisms).
Inspired by the statement of Theorem \ref{MainTheorem},
we ask the following birational analogue of Murphy's law:

\begin{question}
  Do all birational types of  MRC quotients of varieties defined over
  $\overline{\mathbb{Q}}$ occur for irreducible components
  of $\mathrm{Hilb}(\AA^n)$ over $\overline{\mathbb{Q}}$?
\end{question}

\vskip 8pt

The path of the proof of Theorem \ref{MainTheorem} can be followed for any well-behaved{\footnote{The
    generic  point of an irreducible component of a Hilbert scheme of projective
    space should represent the generic object of the moduli space of polarized varieties (and
    not lie on any quadrics).}}
moduli space of polarized varieties (instead of the moduli space of curves).
The result is then the existence, for sufficiently large $n$, of
irreducible components of $\mathrm{Hilb}(\AA^n)$  which dominate these moduli spaces.
The investigation can be thought of as realizing the vision of William Blake in
his poem {\em Auguries of Innocence} (1803):

\vspace{8pt}
{\em \hspace{40pt} To see a world in a grain of sand}\\

\vspace{-20pt}
{\em \hspace{40pt} And a heaven in a wild flower,}\\

\vspace{-20pt}
{\em \hspace{40pt} Hold infinity in the palm of your hand}\\

\vspace{-20pt}
{\em \hspace{40pt} And eternity in an hour.\em}

\section*{Acknowledgments}
Our work started at the
{\em Hilbert schemes of points} conference
at Humboldt University in Berlin in September 2023. The main idea for the construction began in a conversation at
Café Bravo in Berlin-Mitte with Joachim Jelisiejew (who helped us significantly).
Further progress was made in  Les Diablerets in January 2024 at the
{\em Workshop on the enumerative geometry of the Hilbert scheme of points}. We thank
Michele Graffeo, Paolo Lella, Sergej Monavari, and Ritvik Ramkumar for related discussions about Hilbert schemes.

\section*{Funding}
Farkas was supported by the Berlin Mathematics Research Center $\mathrm{MATH+}$ and by the ERC Advanced Grant SYZYGY (834172). Pandharipande was supported by SNF-200020-182181, SNF-200020-219369, ERC-2017-AdG-786580MACI, and by SwissMAP.
Sammartano  was  supported by the grants PRIN 2020355B8Y
{\em Square-free Gr\"obner degenerations, special varieties and related topics},
PRIN 2022K48YYP \emph{Unirationality, Hilbert schemes, and singularities},
and  INdAM – GNSAGA  CUP E55F22000270001.

\end{document}